\documentclass[reqno,12pt,letterpaper]{amsart}
\usepackage[proof]{sdmacros}

\title{Around quantum ergodicity}
\author{Semyon Dyatlov}
\email{dyatlov@math.mit.edu}
\address{Department of Mathematics, Massachusetts Institute of Technology, Cambridge, MA 02139}

\begin{document}

\begin{abstract}
We discuss Shnirelman's Quantum Ergodicity Theorem, giving an outline of a proof
and an overview of some of the recent developments in mathematical Quantum Chaos.

\smallskip

\noindent\textsc{R\'esum\'e.} Nous discuterons le Th\'eor\`eme Ergodicit\'e
Quantique de Shnirelman. Nous donnons l'esquisse d'une preuve
et un aper\c cu des r\'esultats plus r\'ecents dans le domaine math\'ematique
du Chaos Quantique. 
\end{abstract}

\maketitle

\addtocounter{section}{1}
\addcontentsline{toc}{section}{1. Introduction}

Let $(M,g)$ be a compact smooth Riemannian manifold without boundary. Consider an orthonormal basis
of eigenfunctions of the Laplace--Beltrami operator $\Delta_g$:
$$
-\Delta_g u_j=\lambda_j^2 u_j,\quad
\|u_j\|_{L^2}=1,
$$
where the sequence of square roots of eigenvalues $\lambda_j\geq 0$, counted with multiplicity, goes to infinity as $j\to\infty$.

We say that some subsequence $\{u_{j_k}\}$ \emph{equidistributes in physical space}
if the probability measures $|u_{j_k}|^2\,d\vol_g$ converge weakly
to the normalized volume measure, namely
for all $a\in C^\infty(M)$
\begin{equation}
  \label{e:equid-physical}
\int_M a(x)|u_{j_k}(x)|^2\,d\vol_g(x)\to {1\over \vol_g(M)}\int_M a(x)\,d\vol_g(x)\quad\text{as}\quad k\to\infty.
\end{equation}
In~\cite{Shnirelman1} Shnirelman announced the following remarkable theorem
(or rather, its more general version similar to Theorem~\ref{t:QE-phase} below)
which is one of the foundational results in the field of Quantum Chaos:
\begin{theo}[Shnirelman's Theorem/Quantum Ergodicity]
  \label{t:QE-physical}
Assume that the geodesic flow on~$M$ is ergodic with respect to the Liouville measure
(see~\eqref{e:ergodic-def} below). Then
there exists a density~1 subsequence $\{u_{j_k}\}$ which equidistributes in physical space.
\end{theo}
Here `density 1' means that (once again counting eigenvalues with multiplicity)
$$
{\#\{k\mid \lambda_{j_k}\leq R\}\over \#\{j\mid \lambda_j\leq R\}}\to 1\quad\text{as}\quad R\to \infty.
$$
Theorem~\ref{t:QE-physical} is striking for two reasons. First of all, the equidistribution property
gives us strong information on the distribution of mass of eigenfunctions at high frequency,
and it has a natural physical interpretation: if we think of $u_{j_k}$ as pure states
of a quantum particle on~$M$, then equidistribution means that in the high energy limit
the probability of finding the particle in a macroscopic set converges to the volume of that set.
Secondly, the assumption made on the chaotic behavior of the geodesic flow is very weak,
and there are plenty of examples of manifolds with ergodic geodesic flows, such as manifolds
of negative sectional curvature. There are also cases where Theorem~\ref{t:QE-physical}
holds but the associated flow is not ergodic, such as rational polygons, see Marklof--Rudnick~\cite{Marklof-Rudnick}. These examples however do not satisfy the stronger Theorem~\ref{t:QE-phase} stated below.

Following the announcement~\cite{Shnirelman1}, more details of the proof were provided in
the lecture notes~\cite{Shnirelman2}.
A detailed proof was given by Shnirelman later in the addendum to~\cite{Lazutkin-Book}.
In between~\cite{Shnirelman2} (which was not accessible in the West
at the time) and~\cite{Lazutkin-Book} two other proofs were produced
by Zelditch~\cite{Zelditch-QE} and Colin de Verdi\`ere~\cite{CdV-QE}. The general principles used
in all of these proofs are similar, and a version of the proof is sketched in~\S\ref{s:proof-outline} below.

\begin{figure}
\includegraphics[height=7.75cm]{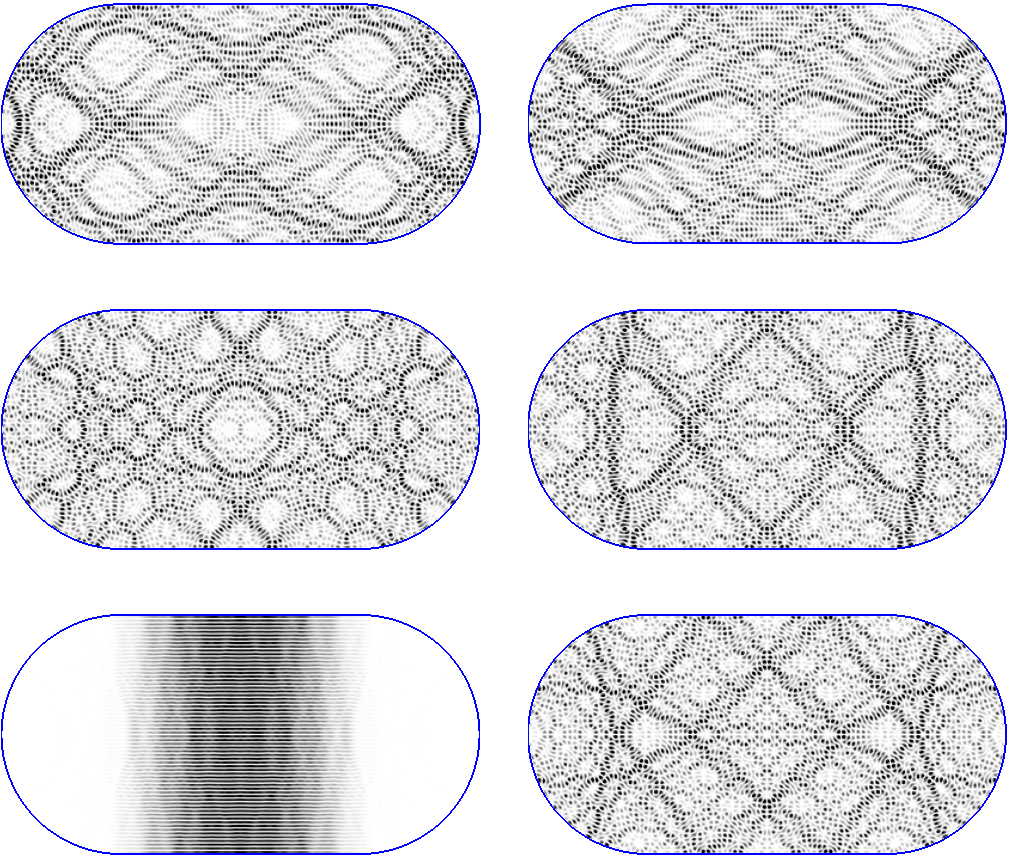}
\qquad
\includegraphics[height=7.75cm]{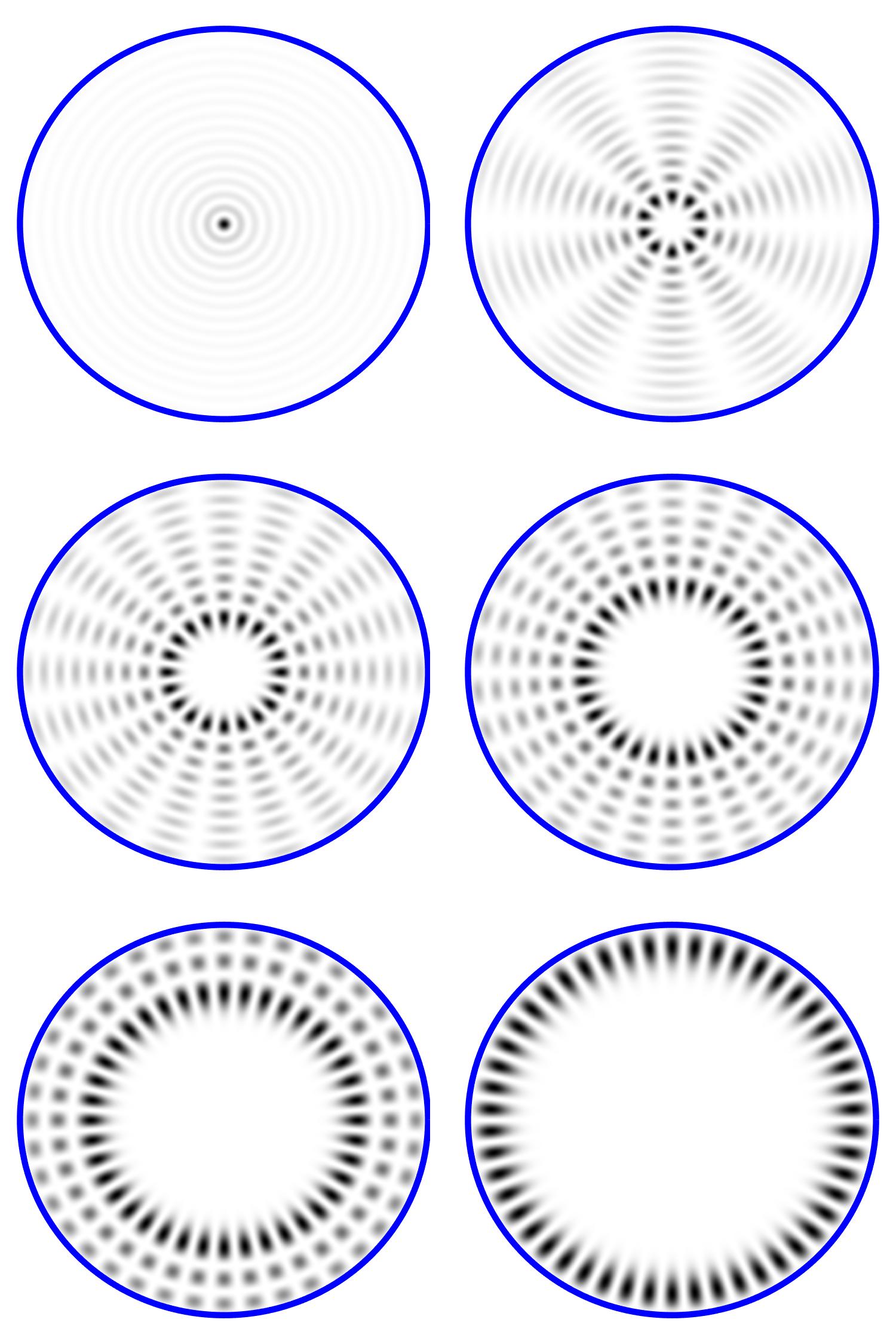}
\caption{Left: several high frequency Dirichlet eigenfunctions of a Bunimovich stadium.
The billiard flow is ergodic and most
eigenfunctions equidistribute as follows from Quantum Ergodicity. The picture is courtesy of Alex Barnett, see~\cite{Barnett-Billiard} and~\cite{Barnett-Hassell} for a description of the method
used and for a numerical investigation of Quantum Ergodicity. For a generic
stadium there is a sequence which does not equidistribute in phase space as proved by Hassell~\cite{Hassell-Bunimovich},
see~\S\ref{s:QUE} below. Right: eigenfunctions of a disk, where the billiard flow is not ergodic
and Quantum Ergodicity fails.
An interesting question is what happens for \emph{mixed systems} which have both ergodic
and non-ergodic regions; see
Schubert~\cite{Schubert-PhD},
Galkowski~\cite{Galkowski-mushrooms}, Rivi\`ere~\cite{Riviere-mushrooms},
and Gomes~\cite{Gomes-mushrooms}.
}
\label{f:barnett}
\end{figure}

A semiclassical version of quantum ergodicity, applying to a larger
class of operators than $-\Delta_g$, was proved
by Helffer--Martinez--Robert~\cite{Helffer-Martinez-Robert}.
In the case of (Dirichlet or Neumann) eigenfunctions
the analog of Theorem~\ref{t:QE-physical} was proved by G\'erard--Leichtnam~\cite{Gerard-Leichtnam} for convex domains in $\mathbb R^n$ with $W^{2,\infty}$ boundaries
(such as the Bunimovich stadium in Figure~\ref{f:barnett})
and by Zelditch--Zworski~\cite{Zelditch-Zworski} for compact
Riemannian manifolds with piecewise $C^\infty$ boundaries. In that
setting the geodesic flow is replaced by the billiard ball flow, which is defined almost
everywhere with respect to the Liouville measure. See Figure~\ref{f:barnett} for a numerical
illustration of quantum ergodicity in this setting. For an overview of various related results,
see~\S\ref{s:overview}.

\subsection{Semiclassical quantization and quantum ergodicity in phase space}

The proofs of Theorem~\ref{t:QE-physical} in fact give a stronger statement, Theorem~\ref{t:QE-phase} below. To state it,
we introduce the notion of \emph{semiclassical quantization}.
To each smooth compactly supported function $a(x,\xi)\in \CIc(T^*M)$ on the cotangent bundle $T^*M$
(called a \emph{classical observable}) we associate the family of operators
(called the corresponding \emph{quantum observable})
$$
\Op_h(a)=a^{\mathrm w}(x,\textstyle{h\over i}\partial_x):L^2(M)\to L^2(M),
$$
depending on the \emph{semiclassical parameter} $0<h\ll 1$. The notation
$a^{\mathrm w}(x,\textstyle{h\over i}\partial_x)$ (where `w' stands for `Weyl') is formal: the operators
$x$ and~$\partial_x$ do not commute, and on a manifold the operators~$\partial_x$
are coordinate dependent. One way to define $\Op_h(a)$ is to start with the Weyl quantization formula for the case $M=\mathbb R^n$
\begin{equation}
  \label{e:weyl-quantization}
\Op_h(a)f(x)=(2\pi h)^{-n}\int_{\mathbb R^n}e^{{i\over h}(x-y)\xi}a(\textstyle{x+y\over 2},\xi)f(y)\,dyd\xi
\end{equation}
and use coordinate charts to piece together a (non-canonical) quantization procedure on a general manifold.

The quantization procedure has many useful properties, in particular the following product and adjoint formulas
for all $a,b\in\CIc(T^*M)$:
\begin{align}
  \label{e:product}
\Op_h(a)\Op_h(b)&=\Op_h(ab)+\mathcal O(h)_{L^2\to L^2};\\
  \label{e:adjoint}
\Op_h(a)^*&=\Op_h(\overline a)+\mathcal O(h)_{L^2\to L^2}.
\end{align}
Moreover, if $a\in \CIc(T^*M)$ then $\|\Op_h(a)\|_{L^2\to L^2}$ is bounded uniformly as $h\to 0$.
We refer the reader to the book of Zworski~\cite{Zworski-Book} for details.

We will choose the semiclassical parameter depending on the eigenvalue $\lambda_j^2$ as follows:
$$
h_j:={1\over\lambda_j}.
$$
This choice is motivated by the observation that the eigenfunction $u_j$ is expected to oscillate
at frequency $\sim \lambda_j$ and the normalized differential operators $h_j\partial_{x_\ell}$
then roughly preserve the magnitude of $u_j$. Under this choice of $h_j$,
the eigenvalue equation $(-\Delta_g-\lambda_j^2)u_j=0$ becomes
\begin{equation}
   \label{e:Laplace}
(-h_j^2\Delta_g-1)u_j=0.
\end{equation}
One can define $\Op_h(a)$ for observables~$a$ with controlled growth as $\xi\to\infty$
rather than just compactly supported ones (see~\cite{Zworski-Book}). In particular, if $a=a(x)$ is a function of~$x$ only, then
$\Op_h(a)$ is the multiplication operator by~$a$, which means that the left-hand side of~\eqref{e:equid-physical}
can be written as follows:
$$
\int_M a(x)|u_j(x)|^2\,d\vol_g(x)=\langle \Op_{h_j}(a)u_j,u_j\rangle\quad\text{for all}\quad
a\in C^\infty(M)
$$
where $\langle \bullet,\bullet\rangle$ denotes the inner product on $L^2(M,d\vol_g)$.

For general classical observables $a(x,\xi)$, the expression $\langle\Op_h(a)u,u\rangle$ can be interpreted as the average value of $a$ for a quantum particle with wave function~$u$;
here $x$ denotes the position variables and $\xi$ the momentum variables.
This suggests the following generalization of~\eqref{e:equid-physical}: we say that
$\{u_{j_k}\}$ \emph{equidistributes in phase space} if
for all $a\in \CIc(T^*M)$ we have
\begin{equation}
  \label{e:equid-phase}
\langle \Op_{h_{j_k}}(a)u_{j_k},u_{j_k}\rangle\to \int_{S^*M} a(x,\xi)\,d\mu_L(x,\xi)\quad\text{as}\quad
k\to\infty.
\end{equation}
Here the Liouville measure~$\mu_L$ on the cosphere bundle $S^*M=\{(x,\xi)\in T^*M\colon |\xi|_g=1\}$ is defined by
$$
d\mu_L(x,\xi) = c\,d\vol_g(x) dS(\xi)
$$
where the densities $d\vol_g$ on $M$ and $dS$ on the fibers of~$S^*M$ are induced by the metric~$g$ and the constant $c>0$ is chosen
so that $\mu_L$ be a probability measure.

The restriction to $S^*M$ in~\eqref{e:equid-phase}
comes from the fact that eigenfunctions `live' on the
cosphere bundle, more precisely
\begin{equation}
  \label{e:elliptic}
a\in\CIc(T^*M),\quad a|_{S^*M}=0\quad\Longrightarrow\quad
\|\Op_{h_j}(a)u_j\|_{L^2}=\mathcal O(h_j).
\end{equation}
To see~\eqref{e:elliptic} we write $-h^2\Delta_g-1=\Op_h(|\xi|_g^2-1)+\mathcal O(h)$,
which together with the product formula~\eqref{e:product}
(or rather, its version for symbols that are not compactly supported)
gives $\Op_h(a)=\Op_h(b)(-h^2\Delta_g-1)+\mathcal O(h)_{L^2\to L^2}$
where $b:=(|\xi|_g^2-1)^{-1}a\in \CIc(T^*M)$. Applying this to $u_j$
and using the eigenvalue equation~\eqref{e:Laplace} we get~\eqref{e:elliptic}.

The semiclassical version of Theorem~\ref{t:QE-physical} is now given by
\begin{theo}[Quantum Ergodicity in phase space]
  \label{t:QE-phase}
Assume that the geodesic flow on~$M$ is ergodic with respect to the Liouville measure. Then
there exists a density~1 subsequence $\{u_{j_k}\}$ which equidistributes in phase space in the sense of~\eqref{e:equid-phase}.
\end{theo}
Here the geodesic flow is considered as a flow on the cosphere bundle, denoted by
$$
\varphi_t:S^*M\to S^*M,
$$
and ergodicity is defined as follows: for any Borel set $U\subset S^*M$
\begin{equation}
  \label{e:ergodic-def}
\varphi_t(U)=U\quad\text{for all }t\quad\Longrightarrow\quad 
\mu_L(U)=0\quad\text{or}\quad\mu_L(U)=1.
\end{equation}
We note that equidistribution in phase space is a stronger property than equidistribution in the physical
space. A basic example is when $M=\mathbb R/2\pi\mathbb Z$ is a circle, then
the sequence of eigenfunctions $u_j=e^{ijx}$ equidistributes in the physical space
but it does not equidistribute in the phase space: instead it is localized
on the `positive half' of the cosphere bundle, given by $\{\xi=1\}$.

\section{Outline of the proof of quantum ergodicity}
  \label{s:proof-outline}

We now give an outline of the proof of Theorem~\ref{t:QE-phase} to illustrate
the main ideas used. Our presentation roughly follows the book of Zworski~\cite[Chapter~15]{Zworski-Book}
and lecture notes by the author~\cite{qe-notes}
and we refer to these sources for the details omitted here.
This strategy of the proof is due to Zelditch~\cite{Zelditch-C*} in the abstract setting
of $C^*$-algebras.

\subsection{Reduction to an averaged statement}

We reduce the `density 1' type statement of Theorem~\ref{t:QE-phase}
to an estimate averaged over eigenfunctions, Theorem~\ref{t:QE-averaged} below.
To do this we use the \emph{Weyl law}, which gives the asymptotic growth of eigenvalues:
\begin{equation}
  \label{e:Weyl-law}
\#\{j\mid \lambda_j\leq R\}={\omega_n\over (2\pi)^n}\vol_g(M)R^n+o(R^n)\quad\text{as}\quad R\to\infty
\end{equation}
where $n=\dim M$ and $\omega_n>0$ is the volume of the unit ball in~$\mathbb R^n$.

The Weyl law~\eqref{e:Weyl-law} can be proved using two properties
of semiclassical quantization:
\begin{itemize}
\item \emph{Functional calculus}: If $\chi\in \CIc(\mathbb R)$ then
$\chi(-h^2\Delta_g)=\Op_h(a_\chi)$ is the quantization of a symbol
$a_\chi(x,\xi;h)\in \CIc(T^*M)$ which has a full expansion in powers of~$h$
and $a_\chi(x,\xi;h)=\chi(|\xi|_g^2)+\mathcal O(h)$ as $h\to 0$.
See~\cite[Theorem~14.9]{Zworski-Book}.
\item \emph{Trace formula}: If $a\in \CIc(T^*M)$ then the operator
$\Op_h(a):L^2(M)\to L^2(M)$ is trace class and
\begin{equation}
  \label{e:trace-formula}
\tr \Op_h(a)=(2\pi h)^{-n}\bigg( \int_{T^*M} a(x,\xi)\,d\xi dx+\mathcal O(h)\bigg)\quad\text{as}\quad h\to 0
\end{equation}
where $d\xi dx$ is the (canonically defined) symplectic volume form on $T^*M$.
Note that in the case of Weyl quantization on $\mathbb R^n$ defined by~\eqref{e:weyl-quantization} it is easy to see that the trace formula is exact, using that
the trace of an operator is the integral of its Schwartz kernel on the diagonal.
\end{itemize}
Combining these two statements gives for any $\chi\in \CIc(\mathbb R)$
\begin{equation}
  \label{e:Weyl-trace}
\sum_j \chi\big((h\lambda_j)^2\big)=\tr\chi(-h^2\Delta_g)=(2\pi h)^{-n}\bigg(\int_{T^*M} \chi(|\xi|_g^2)\,d\xi dx+\mathcal O(h)\bigg)
\end{equation}
and the Weyl law follows by taking $h:=R^{-1}$, approximating the indicator function $\mathbf 1_{[0,1]}$ above and below
by smooth functions $\chi$, and using the monotonicity of the left-hand side of~\eqref{e:Weyl-trace} in~$\chi$.

We now state the integrated form of Quantum Ergodicity:
\begin{theo}[Integrated Quantum Ergodicity]
  \label{t:QE-averaged}
Assume that the geodesic flow on~$M$ is ergodic with respect to the Liouville measure.
Take arbitrary $a\in \CIc(T^*M)$ and put
\begin{equation}
  \label{e:V-j-def}
L_a:=\int_{S^*M} a\,d\mu_L,\qquad
V_j(a):=\langle\Op_{h_j}(a)u_j,u_j\rangle
\end{equation}
where $\mu_L$ is the Liouville measure on $S^*M$. Then
$$
R^{-n}\sum_{\lambda_j\in [R,2R]}|V_j(a)-L_a|^2\to 0\quad\text{as}\quad R\to \infty.
$$ 
\end{theo}
Theorem~\ref{t:QE-averaged} states that if we restrict $\lambda_j$
to a spectral window $[R,2R]$, where $R$ is large, then $|V_j(a)-L_a|$ is small
on average. By the Chebyshev inequality, we then see that there exists
$\varepsilon(R)$ which goes to zero as $R\to\infty$ such that
$|V_j(a)-L(a)|\leq \varepsilon(R)$ for all $\lambda_j\in [R,2R]$
except an $\varepsilon(R)$ proportion of these. Taking $R=2^n$, $n\to\infty$,
we get a density 1 sequence of $u_j$ which equidistribute in the phase space
for a given classical observable~$a$. Using a diagonal argument one
can construct a density 1 sequence which equidistributes with respect to every
observable, thus giving Theorem~\ref{t:QE-phase}.

In the remainder of this section we sketch a proof of Theorem~\ref{t:QE-averaged}.
We restrict ourselves to the special case when $L_a=0$:
\begin{equation}
  \label{e:QE-ave}
L_a=0\quad\Longrightarrow\quad R^{-n}\sum_{\lambda_j\in [R,2R]}|V_j(a)|^2\to 0\quad\text{as}\quad R\to\infty.
\end{equation}
The general case follows by applying~\eqref{e:QE-ave} to the shifted observable
$a-L_a$ and noting that
$V_j(a-L_a)=V_j(a)-L_a$ since $\Op_h(1)$ is the identity operator.
(Here $a-L_a$ is not compactly supported but this does not make a difference in the proof.)

\subsection{Replacing by ergodic averages}

We will next use the semiclassical Schr\"odinger propagator, which is the unitary family of operators
$$
U(t)=U(t;h):=\exp(ith\Delta_g/2):L^2(M)\to L^2(M).
$$
For any $a\in \CIc(T^*M)$ we have
\begin{equation}
  \label{e:v-j-prop}
V_j(a)=\langle \Op_{h_j}(a)u_j,u_j\rangle=\langle U(-t;h_j)\Op_{h_j}(a)U(t;h_j)u_j,u_j\rangle.
\end{equation}
Here we use that $u_j$ is an eigenfunction of the Laplacian and thus also
of~$U(t)$, more precisely $U(t;h_j)u_j=e^{-it\lambda_j/2}u_j$. In fact, this is the most important
place where one uses the fact that $u_j$'s are eigenfunctions.

The conjugated operator $U(-t;h)\Op_h(a)U(t;h)$ is described by \emph{Egorov's Theorem}:
\begin{equation}
  \label{e:egorov}
U(-t;h)\Op_h(a)U(t;h)=\Op_h(a\circ\varphi_t)+\mathcal O_t(h)_{L^2\to L^2}
\end{equation}
where $\varphi_t$, defined as the Hamiltonian flow of $|\xi|^2_g/2$, is an extension of the geodesic flow from $S^*M$
to $T^*M$ and the constant in the remainder depends on $t$ but not on~$j$.
The statement~\eqref{e:egorov} is what relates classical dynamics (the geodesic flow)
to quantum dynamics (the Schr\"odinger propagator) in the proof. See~\cite[Theorem~15.2]{Zworski-Book} for details.

Combining~\eqref{e:v-j-prop} and~\eqref{e:egorov} we see that
the quantum observables $V_j(a)$ do not change much if we replace $a$
by its pullback by the geodesic flow:
$$
V_j(a)=V_j(a\circ\varphi_t)+\mathcal O_t(h_j)\quad\text{for all}\quad t\in\mathbb R.
$$
Since $V_j$ depends linearly on $a$, same is then true for \emph{ergodic averages}:
\begin{equation}
  \label{e:V-j-avg}
V_j(a)=V_j(\langle a\rangle_T)+\mathcal O_T(h_j)\quad\text{for all}\quad T>0
\end{equation}
where the ergodic average $\langle a\rangle_T$ is defined by
$$
\langle a\rangle_T={1\over T}\int_0^T a\circ\varphi_t\,dt\ \in\ \CIc(T^*M).
$$
Recall that we assumed that the flow $\varphi_t$ is ergodic on $S^*M$ with respect
to the Liouville measure $\mu_L$ and that $\int_{S^*M}a\,d\mu_L=0$.
Then by the von Neumann Ergodic Theorem we have
\begin{equation}
  \label{e:von-neumann}
\|\langle a\rangle_T\|_{L^2(S^*M;d\mu_L)}\to 0\quad\text{as}\quad T\to \infty.
\end{equation}
It remains to convert this bound on $\langle a\rangle_T$
to a bound on $V_j(\langle a\rangle_T)$.

\subsection{Local Weyl bound}

We now bound the right-hand side of~\eqref{e:QE-ave}:
\begin{lemm}
  \label{l:HS-bound}
There exists a constant $C$ such that for all $a\in\CIc(T^*M)$ and all $R>1$
\begin{equation}
\label{e:HS-bound}
R^{-n}\sum_{\lambda_j\in [R,2R]}|V_j(a)|^2\leq C\|a\|_{L^2(S^*M;d\mu_L)}^2+\mathcal O_{a}(R^{-1}).
\end{equation}
\end{lemm}
Roughly speaking, Lemma~\ref{l:HS-bound} says that if $a(x,\xi)$ is small after averaging over the points $(x,\xi)$ then $V_j(a)$ is small after averaging over the eigenvalues $\lambda_j$.

To show Lemma~\ref{l:HS-bound},
we first bound the sum in~\eqref{e:HS-bound} by a smoothened out version.
Fix a nonnegative cutoff function $\chi\in \CIc((0,\infty))$ such that $\chi=1$ on $[1,4]$.
Since
$$
|V_j(a)|=|\langle \Op_{h_j}(a)u_j,u_j\rangle|\leq \|\Op_{h_j}(a)u_j\|_{L^2},
$$
we have
\begin{equation}
  \label{e:long-formula}
\begin{aligned}
R^{-n}\sum_{\lambda_j\in [R,2R]}|V_j(a)|^2&\leq R^{-n}\sum_j\chi\Big({\lambda_j^2\over R^2}\Big)
\|\Op_{h_j}(a)u_j\|_{L^2}^2\\
&=R^{-n}\sum_j \chi\Big({\lambda_j^2\over R^2}\Big)\langle \Op_{h_j}(|a|^2)u_j,u_j\rangle+\mathcal O(R^{-1}).
\end{aligned}
\end{equation}
Here in the last line we use that $\Op_h(a)^*\Op_h(a)=\Op_h(|a|^2)+\mathcal O(h)_{L^2\to L^2}$
as follows from the algebraic properties~\eqref{e:product}, \eqref{e:adjoint};
we also use the Weyl law~\eqref{e:Weyl-law}.

The first term on the right-hand side of~\eqref{e:long-formula} looks very much like a trace
except that $\Op_{h_j}(|a|^2)$ depends on~$j$. To remove this dependence, we rescale
the semiclassical parameter. Put
$$
h:={1\over R},\quad
\tau_j:={h\over h_j}={\lambda_j\over R},\quad
\tau_j^2\in\supp\chi.
$$
For any $b\in\CIc(T^*M)$ we have (here we abuse the formal notation $b^{\mathrm w}(x,{h\over i}\partial_x)$;
see~\cite[Exercise~E.5]{DZ-Book} for more details)
$$
\Op_{h}(b)=b^{\mathrm w}(x,\textstyle{h\over i}\partial_x)=
b^{\mathrm w}(x,\tau_j\textstyle{h_j\over i}\partial_x)=\Op_{h_j}(\Lambda_{\tau_j} b)
\quad\text{where}\quad \Lambda_\tau b(x,\xi):=b(x,\tau\xi).
$$
This implies that $V_j(\Lambda_{\tau_j}b)=\langle \Op_h(b)u_j,u_j\rangle$.
Now, fix $b$ such that $\Lambda_\tau b=|a|^2$ on $S^*M$ for all $\tau>0$
such that $\tau^2\in\supp\chi$; for example, one can put
$b(x,r\theta):=\psi(r)|a(x,\theta)|^2$ for all $(x,\theta)\in S^*M$,
$r>0$, and an appropriate choice of the cutoff $\psi\in\CIc((0,\infty))$.
By~\eqref{e:elliptic} we have $\|\Op_{h_j}(|a|^2-\Lambda_{\tau_j}b)u_j\|_{L^2}=\mathcal O(R^{-1})$
and thus
$$
V_j(|a|^2)=V_j(\Lambda_{\tau_j}b)+\mathcal O(R^{-1})=\langle\Op_h(b)u_j,u_j\rangle+\mathcal O(R^{-1})
$$
where on the right-hand side the semiclassical parameter $h:=R^{-1}$
no longer depends on~$j$. Putting this together
with~\eqref{e:long-formula} and the Weyl law~\eqref{e:Weyl-law}, we get
$$
\begin{aligned}
R^{-n}\sum_{\lambda_j\in [R,2R]}|V_j(a)|^2
&\leq R^{-n}\sum_j \chi\Big({\lambda_j^2\over R^2}\Big)
\langle\Op_h(b)u_j,u_j\rangle+\mathcal O(R^{-1})\\
&=R^{-n}\tr\big(\chi(-h^2\Delta_g)\Op_h(b)\big)+\mathcal O(R^{-1}).
\end{aligned}
$$
By the functional calculus and the product formula~\eqref{e:product}
we have $\chi(-h^2\Delta_g)\Op_h(b)=\Op_h(\chi(|\xi|_g^2)b)+\mathcal O(h)$.
By the trace formula~\eqref{e:trace-formula} we then have
$$
\begin{aligned}
R^{-n}\sum_{\lambda_j\in [R,2R]}|V_j(a)|^2&\leq (2\pi)^{-n}\int_{T^*M} \chi(|\xi|_g^2)b(x,\xi)\,d\xi dx+\mathcal O_a(R^{-1})\\&\leq C \|a\|_{L^2(S^*M;d\mu_L)}^2+\mathcal O_a(R^{-1})
\end{aligned}
$$
which finishes the proof of~\eqref{e:HS-bound}.

\subsection{End of the proof}

We are now ready to finish the proof of~\eqref{e:QE-ave}, and thus of Theorem~\ref{t:QE-phase}.
Take some large $T>0$. We have
$$
\begin{aligned}
R^{-n}\sum_{\lambda_j\in [R,2R]}|V_j(a)|^2 &
= R^{-n}\sum_{\lambda_j \in [R,2R]} |V_j(\langle a\rangle_T)|^2 + \mathcal O_T(R^{-1})
\\
& \leq C\|\langle a\rangle_T\|^2_{L^2(S^*M;d\mu_L)} + \mathcal O_T(R^{-1}).
\end{aligned}
$$
Here in the first line we used~\eqref{e:V-j-avg} and the Weyl law~\eqref{e:Weyl-law};
in the second line we used~\eqref{e:HS-bound}. The constant in $\mathcal O(\bullet)$
depends on~$T$ but the constant~$C$ does not.

Passing to the limit $R\to\infty$, we get
$$
\limsup_{R\to\infty}R^{-n}\sum_{\lambda_j\in [R,2R]}|V_j(a)|^2
\leq C\|\langle a\rangle_T\|^2_{L^2(S^*M;d\mu_L)}.
$$
The left-hand side does not depend on $T$ and the right-hand side
converges to~0 as $T\to\infty$ by~\eqref{e:von-neumann}.
Therefore
$$
\lim_{R\to\infty}R^{-n}\sum_{\lambda_j\in [R,2R]}|V_j(a)|^2=0
$$
which gives~\eqref{e:QE-ave} and finishes the proof.

\section{Overview of literature}
\label{s:overview}

Shnirelman's Quantum Ergodicity Theorem inspired a whole new direction of research
on concentration of eigenfunctions and related objects. Here we give an overview of some of these developments. The list of topics discussed is by no means complete (achieving this would be difficult), and is somewhat skewed towards the author's own research interests.
We refer the reader to the books by Sogge~\cite{Sogge-Hangzhou} and Zelditch~\cite{Zelditch-QE-book},
as well as the review articles by Marklof~\cite{Marklof-QE}, Sarnak~\cite{Sarnak-QE-review},
and Nonnenmacher~\cite{Nonnenmacher-Anatomy},
for more detailed overview of results on quantum ergodicity and related topics.

\subsection{Quantum Unique Ergodicity and semiclassical measures}
\label{s:QUE}

\begin{figure}
\includegraphics[height=6.5cm]{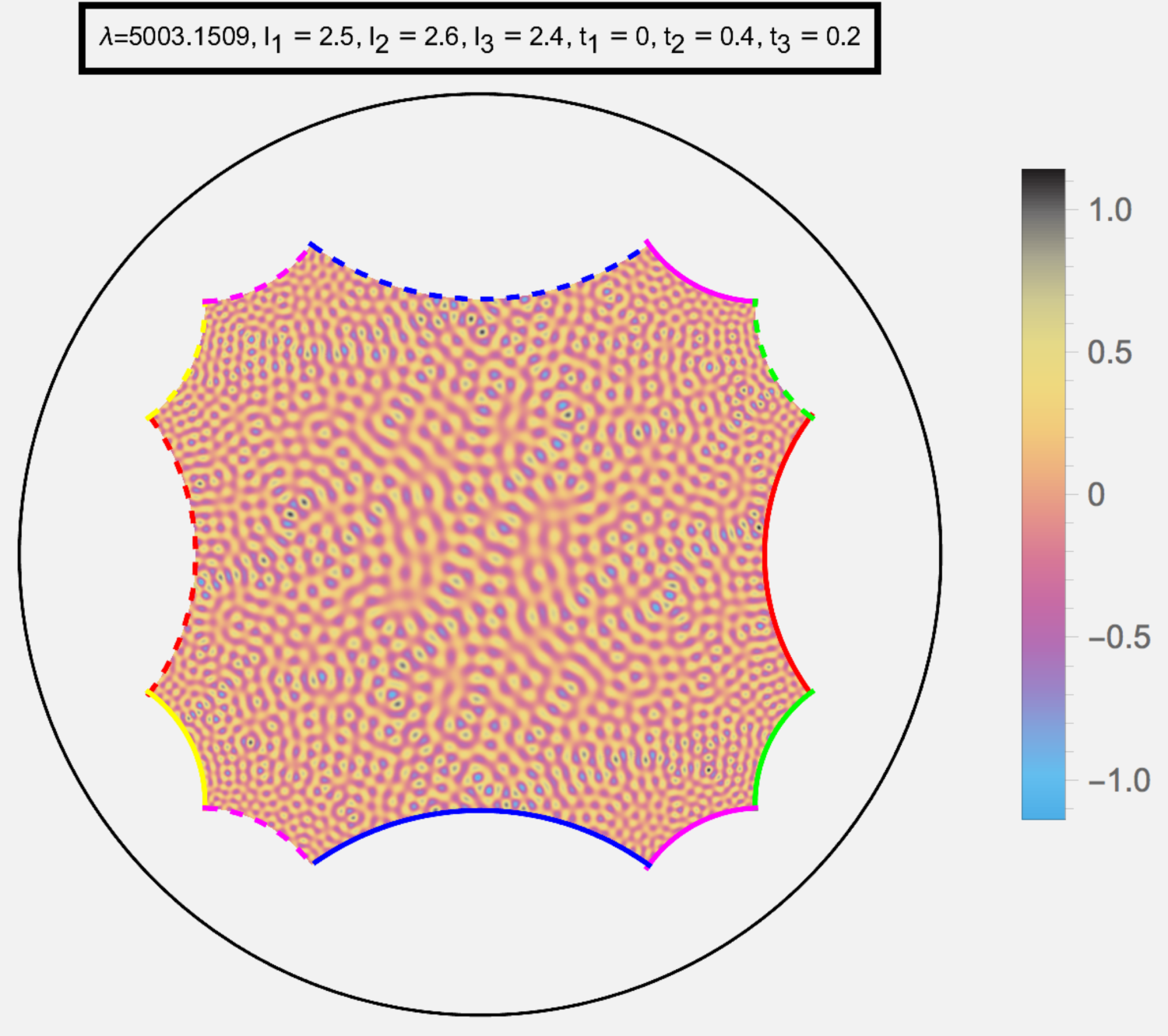}
\quad
\includegraphics[height=6.5cm]{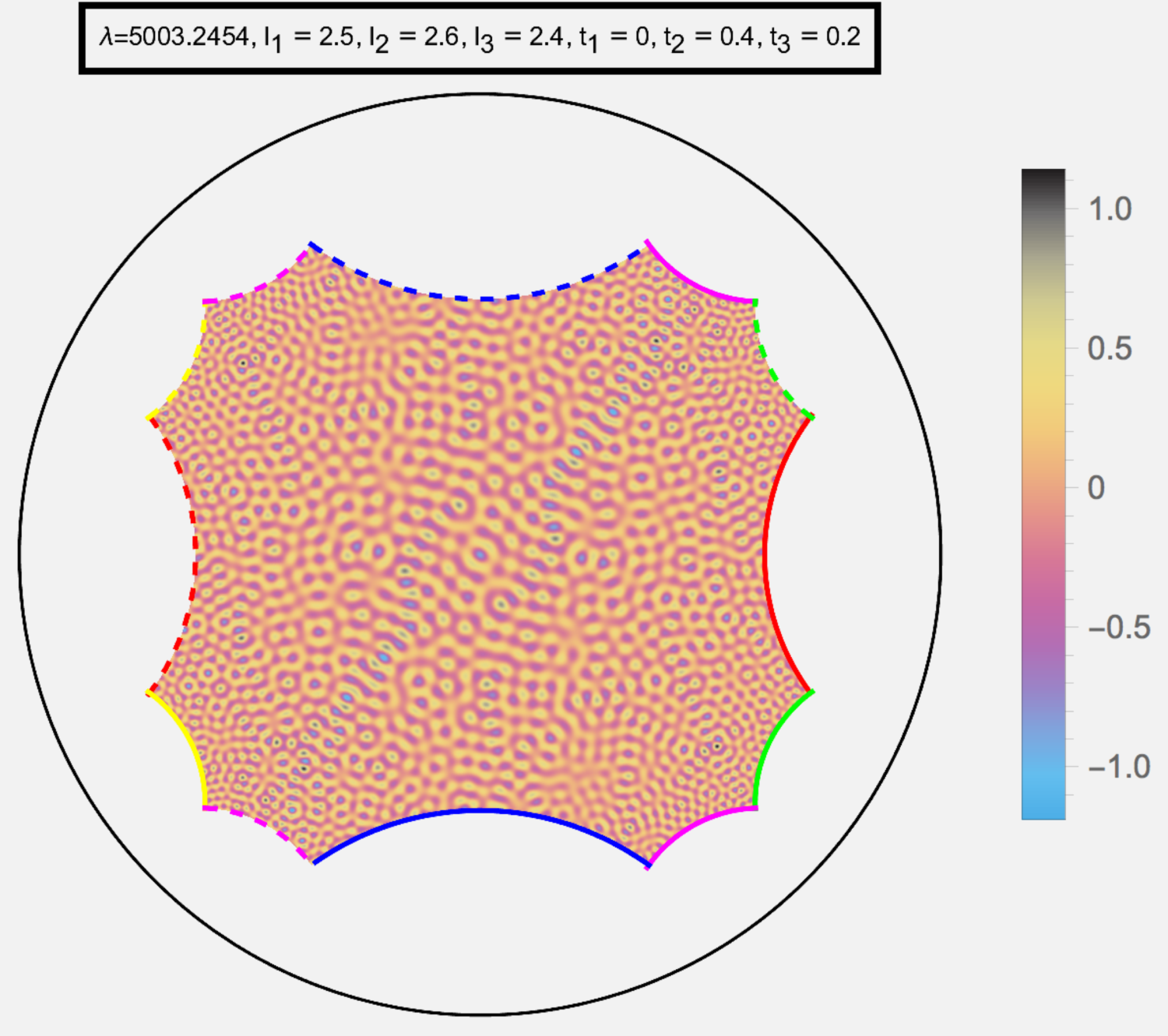}
\caption{Two numerically computed eigenfunctions on a hyperbolic surface $M=\Gamma\backslash\mathbb H^2$,
drawn here on a fundamental domain of the surface in the Poincar\'e disk model of $\mathbb H^2$. While the microscopic features of the two eigenfunctions are different,
on the macroscopic level they both show equidistribution.
Pictures courtesy of Alexander Strohmaier, see Strohmaier--Uski~\cite{Strohmaier-Uski} for more details.}
\label{f:stro}
\end{figure}

Theorem~\ref{t:QE-phase} gives equidistribution in phase space for a density one sequence
of eigenfunctions under a weak chaotic assumption on the geodesic flow. It is natural
to ask if under stronger assumptions one can show that \emph{all} eigenfunctions
equidistribute, i.e. for all $a\in\CIc(T^*M)$
$$
\langle\Op_{h_j}(a)u_j,u_j\rangle\to \int_{S^*M}a\,d\mu_L\quad\text{as}\quad j\to\infty.
$$
This statement, known as \emph{Quantum Unique Ergodicity (QUE)}, was conjectured
by Rudnick--Sarnak~\cite{Rudnick-Sarnak-QUE} for hyperbolic surfaces
(i.e. compact surfaces of constant curvature~$-1$).
Since then there has been much progress (some of which is described below) but the original conjecture
is still very much open. See Figure~\ref{f:stro} for a numerical illustration.
Quantum Unique Ergodicity can also be formulated in the context of \emph{semiclassical measures}, defined as follows:
\begin{defi}
Let $\{u_{j_k}\}$ be a subsequence of eigenfunctions of $-\Delta_g$
and $\mu$ be a probability measure on~$T^*M$.
We say $u_{j_k}$ converges weakly to $\mu$ if
for all $a\in \CIc(T^*M)$ we have
$$
\langle\Op_{h_{j_k}}(a)u_{j_k},u_{j_k}\rangle\to \int_{T^*M}a\,d\mu\quad\text{as}\quad k\to\infty.
$$
We say that $\mu$ is a \emph{semiclassical measure} if it is the weak
limit of some subsequence of eigenfunctions.
\end{defi}
Every semiclassical measure is supported on~$S^*M$
and invariant under the geodesic flow, see for example~\cite[Chapter~5]{Zworski-Book}.
The Quantum Unique Ergodicity conjecture can now be restated as follows:
\emph{the only semiclassical measure is the Liouville measure}.

\subsubsection{The arithmetic case}

A hyperbolic surface can be represented as the quotient $M=\Gamma\backslash\mathbb H^2$
where $\mathbb H^2$ is the hyperbolic plane and $\Gamma$ is a group of isometries.
If we use the upper half-plane model $\mathbb H^2=\{z\in\mathbb C\mid \Im z>0\}$,
with the hyperbolic metric $g={|dz|^2\over |\Im z|^2}$, then orientation preserving isometries
are M\"obius maps
$$
z\mapsto {az+b\over cz+d},\quad
a,b,c,d\in\mathbb R,\quad
ad-bc=1
$$
and the total group of isometries is $\PSL(2,\mathbb R)$, the quotient of $\SL(2,\mathbb R)$
by the group $\{I,-I\}$. Thus compact hyperbolic surfaces are identified with co-compact discrete subgroups $\Gamma\subset\PSL(2,\mathbb R)$.

A special class of hyperbolic surfaces are the \emph{arithmetic} ones, where the group $\Gamma$
has certain number theoretic properties. A particularly important example is the \emph{modular surface} $\PSL(2,\mathbb Z)\backslash\mathbb H^2$, which however is not compact (it has a cusp). For an example of a compact arithmetic hyperbolic surface, see for instance~\cite[\S2]{Marklof-QE}.

Arithmetic surfaces have additional symmetries called \emph{Hecke operators}
(see for example~\cite[\S6]{Marklof-QE}). Those are a family of operators
$M_q:L^2(M)\to L^2(M)$ indexed by positive integers~$q$. The Hecke operators
commute with each other and with the Laplacian, so one can form a basis
of eigenfunctions for the Laplacian which are also eigenfunctions of all $M_q$
(we call this a \emph{Hecke basis}).

Using these additional symmetries, Lindenstrauss~\cite{Lindenstrauss-QUE}
was able to prove the Quantum Unique Ergodicity conjecture for a Hecke basis of eigenfunctions
on compact arithmetic surfaces. Brooks--Lindenstrauss~\cite{Brooks-Lindenstrauss-QUE} extended
this to any joint basis of the Laplacian and a single Hecke operator.
To the author's knowledge it is still an open question
whether Quantum Unique Ergodicity holds for every orthonormal basis
of eigenfunctions of just the Laplacian (which might not be a Hecke basis
if eigenvalues of the Laplacian have multiplicity).

\subsubsection{Ergodic systems where QUE fails}
\label{s:QUE-fails}

Hassell~\cite{Hassell-Bunimovich} (following earlier work by Donnelly~\cite{Donnelly-Bunimovich}) constructed examples of manifolds where the geodesic flow
is ergodic but Quantum Unique Ergodicity fails. These examples
include generic Bunimovich stadia, see Figure~\ref{f:barnett}. A key feature of these is the presence of closed trajectories along which the differential of the geodesic/billiard ball flow grows only polynomially in time; for the Bunimovich stadium these are trajectories bouncing between the top and bottom boundary segments.
We note that hyperbolic surfaces
and, more generally, negatively curved manifolds, do not admit such weakly dispersing trajectories because the geodesic flow has the Anosov property and
the differential of the flow grows exponentially fast in time.

Another family of counterexamples to QUE is for toy models
of \emph{quantum maps}. These are families of matrices of size $N\to\infty$
(where the effective semiclassical parameter is $h:=N^{-1}$) which quantize
ergodic symplectic transformations of an even-dimensional torus.
(Quantum Ergodicity was proved in this setting by Bouzouina--de Bi\`evre~\cite{Bouzouina-deBievre}
and Zelditch~\cite{Zelditch-cat}.)
 In particular,
Faure--Nonnenmacher--de Bi\`evre~\cite{Faure-Nonnenmacher-dB}
showed that for certain two-dimensional quantum cat maps there exists a sequence of eigenstates
which converges to the measure
\begin{equation}
  \label{e:mu-half}
\mu=\textstyle{1\over 2}\mu_L+{1\over 2}\delta_\gamma
\end{equation}
where $\mu_L$ is the volume measure and $\delta_\gamma$ is the delta measure on any
a priori given closed trajectory $\gamma$ of the cat map.
Anantharaman--Nonnenmacher~\cite{Anantharaman-Nonnenmacher-Baker}
considered the Walsh-quantized baker's map and constructed semiclassical
measures supported on fractal sets.

It is possible to introduce analogs of Hecke operators for quantum cat maps.
For two-dimensional quantum cat maps arithmetic Quantum Unique Ergodicity was proved
by Kurlberg--Rudnick~\cite{Kurlberg-Rudnick}. On the other hand 
Kelmer~\cite{Kelmer-cat} constructed
semiclassical measures of Hecke eigenfunctions concentrating on proper submanifolds for certain
higher dimensional quantum cat maps.

\subsubsection{Entropy and support of semiclassical measures}

In the absense of QUE (in non-arithmetic settings), a natural problem is to restrict
as much as possible which flow-invariant probability measures on $S^*M$ can arise
as semiclassical measures. In particular, Colin de Verdi\`ere conjectured in~\cite{CdV-QE}
that for hyperbolic surfaces there cannot be a semiclassical measure supported
on a single closed trajectory. This conjecture was solved in the more general setting
of manifolds with Anosov geodesic flows by Anantharaman~\cite{Anantharaman-Entropy}. 
This was followed by results of Anantharaman--Nonnenmacher~\cite{Anantharaman-Nonnenmacher-Entropy}, Anantharaman--Koch--Nonnenmacher~\cite{Anantharaman-Koch-Nonnenmacher}, Rivi\`ere~\cite{Riviere-Entropy-1,Riviere-Entropy-2}, and Anantharaman--Silberman~\cite{Anantharaman-Silberman}. In particular, \cite{Anantharaman-Nonnenmacher-Entropy} proved a lower bound on the Kolmogorov--Sinai entropy $\mathbf h_{\mathrm{KS}}(\mu)$ of any semiclassical measure $\mu$, which for hyperbolic surfaces is
\begin{equation}
  \label{e:entropy-bound}
\mathbf h_{\mathrm{KS}}(\mu)\geq \textstyle{1\over 2}.
\end{equation}
Here the entropy measures how much of the complexity of the geodesic flow
is captured by $\mu$.
The delta measure $\delta_\gamma$ on a closed geodesic has entropy~0
and the Liouville measure (for hyperbolic surfaces) has entropy~1.
The entropy of the measure of the type~\eqref{e:mu-half} is exactly
$1\over 2$, so from the point of view of the counterexample of~\cite{Faure-Nonnenmacher-dB}
the bound~\eqref{e:entropy-bound} is sharp.

For hyperbolic surfaces, Dyatlov--Jin~\cite{meassupp} showed a different kind of restriction
on semiclassical measures $\mu$: each such measure should have full support,
i.e. $\mu(U)>0$ for any nonempty open set $U\subset S^*M$. (See also the expository article~\cite{Dyatlov-JEDP}.) This rules out the
fractal counterexamples of the kind found in~\cite{Anantharaman-Nonnenmacher-Baker}
(which can have entropy close to~1 and thus are not ruled out by the entropy bound~\eqref{e:entropy-bound}).
There is no contradiction here since the key new ingredient in~\cite{meassupp},
the fractal uncertainty principle of Bourgain--Dyatlov~\cite{fullgap},
does not hold for the Walsh quantization used in~\cite{Anantharaman-Nonnenmacher-Baker}.
On the other hand, linear combinations $c\mu_L+(1-c)\delta_\gamma$ with $0<c<{1\over 2}$
have full support but are ruled out by~\eqref{e:entropy-bound}.
Dyatlov--Jin--Nonnenmacher~\cite{varfup} recently extended the full support property
to general surfaces with Anosov geodesic flows.

\subsection{Quantum ergodicity in other settings}

The natural ideas behind the proof of Theorem~\ref{t:QE-phase} make it very tempting
to try to adapt this proof to other settings beyond compact manifolds. While
each of these settings presents its own unique challenges, and not all the ideas
for the compact case carry well to the more general settings, there have been many
Quantum Ergodicity-style equidistribution results, some of which are briefly reviewed below.

\subsubsection{Surfaces with cusps}

\begin{figure}
\includegraphics[width=15cm]{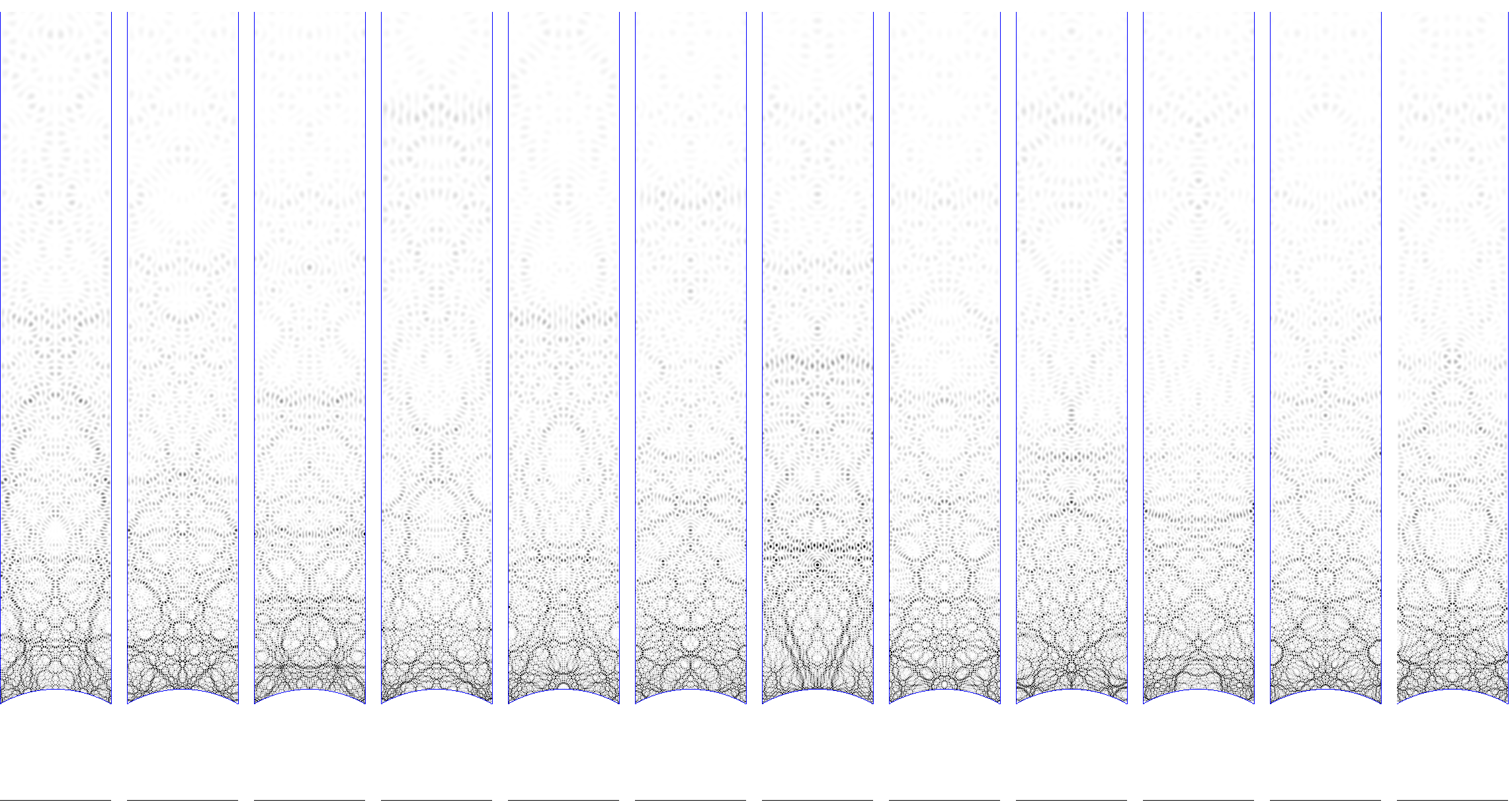}
\caption{A plot of several high-frequency Maass forms on the modular surface,
courtesy of Alex Barnett and made using a code of Holger Then. See \url{https://math.dartmouth.edu/~specgeom/maass.php} and~\cite{Then} for more details.}
\label{f:cusp}
\end{figure}

The first extension of Quantum Ergodicity was to complete noncompact Riemannian surfaces with cusps,
which are infinite ends of the form
$$
[r_0,\infty)_r\times \mathbb S^1_\theta\quad\text{with the metric}\quad
dr^2+e^{-2r}d\theta^2.
$$
An important example is the \emph{modular surface} $\PSL(2,\mathbb Z)\backslash\mathbb H^2$,
which has a fundamental domain of the form $\{|\Re z|\leq {1\over 2},\ |z|\geq 1\}\subset\mathbb H^2$;
here the cusp corresponds to $\Im z\to\infty$. See Figure~\ref{f:cusp}
for a numerical illustration.

The spectrum of the Laplacian on a surface with cusps consists of three parts:
\begin{itemize}
\item The `low' eigenvalues in $[0,{1\over 4})$, which are irrelevant for
high frequency asymptotics featured in Quantum Ergodicity.
\item The continuous spectrum $[{1\over 4},\infty)$.
Assuming for simplicity there is only one cusp, it is parametrized
by \emph{Eisentein functions} $E(x;\lambda)$, $\lambda\geq 0$, which
are certain solutions to the equation
$(-\Delta_g-\lambda^2-{1\over 4})E(x;\lambda)=0$ generalizing one-dimensional plane waves $e^{-ix\lambda}$, $x\in\mathbb R$.
\item The embedded eigenvalues, which are $L^2$ eigenvalues of the Laplacian
in $[{1\over 4},\infty)$. Those are abundant on some surfaces, such as the modular surface,
but generic surfaces with cusps do not have embedded eigenvalues, see Colin de Verdi\`ere~\cite{CdV-pseudolap-1,CdV-pseudolap-2} and Phillips--Sarnak~\cite{Phillips-Sarnak}.
The corresponding eigenfunctions are known as \emph{cusp forms} or \emph{Maass forms}.
\end{itemize}

Zelditch~\cite{Zelditch-QE-cusp} proved Quantum Ergodicity for hyperbolic surfaces with cusps,
of a flavor similar to Theorem~\ref{t:QE-averaged}, featuring both Eisenstein functions and cusp forms. A shorter proof, applying to any surface with cusps which has ergodic geodesic flow, was recently given by Bonthonneau--Zelditch~\cite{Bonthonneau-Zelditch}. In the related setting of eigenfunctions of pseudo-Laplacians Quantum Ergodicity was proved by Studnia~\cite{Studnia-QE}.

In the arithmetic setting of the modular surface, Quantum Unique Ergodicity for
Maass--Hecke forms (that is, Maass forms which are
eigenfunctions of all the Hecke operators) was proved by Soundararajan~\cite{Sound-QUE} following the work of Lindenstrauss~\cite{Lindenstrauss-QUE}. For Eisenstein functions on modular surfaces,
equidistribution was proved by Luo--Sarnak~\cite{Luo-Sarnak-QUE} (in physical space)
and Jakobson~\cite{Jakobson-QUE} (in phase space).

\subsubsection{Manifolds with funnel-type ends}

Another noncompact setting is given by complete Riemannian manifolds $M$ with funnel ends.
The simplest (2-dimensional) version of a funnel end is
\begin{equation}
  \label{e:basic-funnel}
[0,\infty)_r\times \mathbb S^1_\theta\quad\text{with the metric}\quad
dr^2+\cosh^2 r\,d\theta^2.
\end{equation}
In constract with cusp ends, which are very narrow, funnel ends are very wide and in particular
they have infinite volume. For manifolds with funnels, there are no embedded eigenvalues;
the spectrum is purely continuous and parametrized by Eisenstein functions
$E(x;\lambda,\omega)$ where $\lambda$ corresponds to the eigenvalue and
$\omega$ is a point on the conformal infinity of $M$ (for the case of the
basic funnel end~\eqref{e:basic-funnel} we would have $\omega\in\mathbb S^1$).
Alternatively one can consider Euclidean ends, which have the metric
$dr^2+r^2d\theta^2$.

Since $M$ has infinite volume, we can no longer talk about the ergodicity
of the geodesic flow with respect to the Liouville measure.
Instead, one makes assumptions on the set of \emph{trapped geodesics},
which are geodesics which do not escape (forwards or backwards in time)
through the infinite ends of~$M$. This set often has fractal structure.

In the setting of hyperbolic manifolds with funnels,
Guillarmou--Naud~\cite{Guillarmou-Naud-QUE} showed an `equidistribution' statement
for the Eistenstein functions, under the \emph{pressure condition} which here is equivalent
to a certain upper bound on the Hausdorff dimension of the set of trapped geodesics.
This was later extended by Ingremeau~\cite{Ingremeau-QUE} to the setting
where the underlying flow is hyperbolic, assuming again that the pressure condition holds.
Ingremeau later removed the pressure condition in~\cite{Ingremeau-QUE-2}, relying
on the spectral gap proved in~\cite{hgap,fullgap}.
These results are more similar to Quantum Unique Ergodicity than Quantum Ergodicity
in that they are weak convergence statements for the whole family of Eisenstein functions
rather than a density~1 subfamily. The limiting measures depend on $\omega$;
integrating in~$\omega$ one obtains the Liouville measure.

Dyatlov--Guillarmou~\cite{qeefun} proved a Quantum Ergodicity-type statement in this setting
(which applies to most rather than all values of $\lambda,\omega$), replacing the pressure condition with the much weaker assumption that the set of trapped trajectories has Liouville measure~0.

\subsubsection{Restrictions of eigenfunctions}

Coming back to the setting of compact manifolds~$M$, a natural question to ask is the following: for a submanifold $\Sigma\subset M$, do restrictions of a density 1 sequence
of eigenfunctions $u_{j_k}|_{\Sigma}$ converge weakly to a natural measure?
The answer cannot be positive for $\Sigma$ of arbitrary dimension, for example
it is unrealistic to expect equidistribution of restrictions when $\Sigma$ is a point.
Henceforth we restrict to the case when $\Sigma$ is a hypersurface.

It turns out that assuming that $M$ has ergodic geodesic flow is not enough
to ensure equidistribution of restrictions of eigenfunctions. A basic example
is when $\Sigma$ is the fixed point set of some isometric involution $J:M\to M$;
then roughly half of the eigenfunctions $u_j$ are odd with respect to~$J$
and thus have $u_j|_\Sigma=0$.

However, if one makes an additional (generically satisfied)
assumption that $\Sigma$ does not have a reflection symmetry with respect to the geodesic flow,
then $u_{j_k}|_{\Sigma}$ equidistributes for a density one sequence. Here the limiting measure for the equidistribution is naturally defined from the Liouville measure and supported
on the coball bundle $B^*\Sigma=\{(x,\xi)\in T^*\Sigma\colon |\xi|_g\leq 1\}$. This result was proved by Toth--Zelditch~\cite{Toth-Zelditch-1,Toth-Zelditch-2}, with
a different proof of a semiclassical generalization given by Dyatlov--Zworski~\cite{DZ-QE}; see also Christianson--Toth--Zelditch~\cite{Christianson-Toth-Zelditch}.

In case when $M$ is a surface, Quantum Ergodicity for restrictions has applications to counting nodal domains of eigenfunctions, see Jung--Zelditch~\cite{Jung-Zelditch-Nodal-1,Jung-Zelditch-Nodal-2}.

\subsubsection{Large graphs}

Quantum Ergodicity can be also adapted to models of quantum chaos
which are not eigenfunctions of operators on manifolds. In~\S\ref{s:QUE-fails}
we briefly discussed one such setting, quantum maps. Here we briefly discuss
a different setting, \emph{large regular graphs} where quantum ergodicity was proved by Anantharaman--Le Masson~\cite{Anantharaman-LeMasson-Graphs}. We refer the reader to the review of Anantharaman--Sabri~\cite{Anantharaman-Sabri-Review} for more information and further results.

Let $(G_N)$ be a sequence of graphs of size~$N$ which are $k$-regular with
some fixed $k\geq 3$. We identify the set of vertices of $G_N$ with $\{1,\dots,N\}$
and functions on this set with vectors in $\mathbb C^N$.
We will replace the high eigenvalue limit of Theorem~\ref{t:QE-physical}
by the limit $N\to\infty$ and eigenfunctions of the Laplacian by an orthonormal basis $u_j^{(N)}\in \mathbb C^N$, $1\leq j\leq N$, of eigenvectors of the graph Laplacian on $G_N$
(for regular graphs, this is same as eigenvectors of the adjacency matrix of the graph).

The graphs $G_N$ for different values of~$N$ have no relationship with each other,
so it is unclear how to define the quantization of a fixed $N$-independent observable. Instead
one can study the expressions $\langle a^{(N)}u_j^{(N)},u_j^{(N)}\rangle$
for \emph{any} family of functions $a^{(N)}:\{1,\dots,N\}\to\mathbb C$
with $\max |a^{(N)}|\leq 1$. Under two assumptions discussed below,
the paper~\cite{Anantharaman-LeMasson-Graphs} proves the following version
of Integrated Quantum Ergodicity: for any choice of $a^{(N)}$
$$
\lim_{N\to\infty}{1\over N}\sum_{j=1}^N \Big| \langle a^{(N)}u_j^{(N)},u_j^{(N)}\rangle
-\langle a^{(N)}\rangle\Big|^2=0\quad\text{where}\quad
\langle a^{(N)}\rangle:={1\over N}\sum_{\ell=1}^N a^{(N)}(\ell).
$$
The assumptions on the graphs $G_N$ are as follows:
\begin{enumerate}
\item $G_N$ converges to the regular $k$-tree in the local weak sense, that is
the injectivity radius of a random vertex in $G_N$ converges to infinity in distribution;
\item $G_N$ is an expander, that is there is an $N$-independent $\beta>0$ such that
the spectrum of the adjacency matrix of $G_N$ is contained in
$[-k+\beta,k-\beta]\cup \{k\}$.
\end{enumerate}
These assumptions are satisfied for random graphs, as well as for certain deterministic
examples. For random graphs, Bauerschmidt--Huang--Yau~\cite{Bauerschmidt-Huang-Yau}
proved stronger equidistribution statements which are the analog of Quantum Unique Ergodicity in this setting.


\bibliographystyle{alpha}
\bibliography{General,Dyatlov,QC}

\end{document}